\documentclass[
final
]{dmtcs-episciences}

\usepackage{amsmath}
\usepackage{amsthm}
\usepackage{amssymb}
\usepackage{thmtools}
\usepackage{complexity}
\usepackage[shortlabels]{enumitem}
\usepackage[nameinlink,capitalise]{cleveref}
\usepackage{xspace}

\usepackage{tikz}

\newtheorem{theorem}{Theorem}[section]
\newtheorem{lemma}[theorem]{Lemma}
\newtheorem{remark}[theorem]{Remark}
\newtheorem{observation}[theorem]{Observation}
\newtheorem{claim}[theorem]{Claim}
\newtheorem{corollary}[theorem]{Corollary}
\theoremstyle{definition}
\newtheorem{definition}{Definition}

\tikzstyle{td} = [rectangle, draw, fill=yellow!20, text centered, rounded corners]
\tikzstyle{tdarrow} = [->, style = double, double distance=2pt]

\usepackage[utf8]{inputenc}
\usepackage{subfigure}

\usepackage[round]{natbib}

\author[Aboulker et al.]{%
Pierre Aboulker\affiliationmark{1}\thanks{Research supported by french \it{Agence Nationale de la Recherche} under contracts DAGDigDec (JCJC) ANR-21-CE48-0012 and GODASse ANR-24-CE48-4377, and by the group Casino/ENS Chair on Algorithmics and Machine Learning}
\and Nacim Oijid\affiliationmark{2,3}\thanks{Research supported by the Kempe Foundation Grant No. JCSMK24-515 (Sweden).}
\and Robin Petit\affiliationmark{4} \\
\and Mathis Rocton\affiliationmark{5}\thanks{Research supported by the Austrian Science Foundation (FWF, Project 10.55776/Y1329) and the \includegraphics[width=0.4cm]{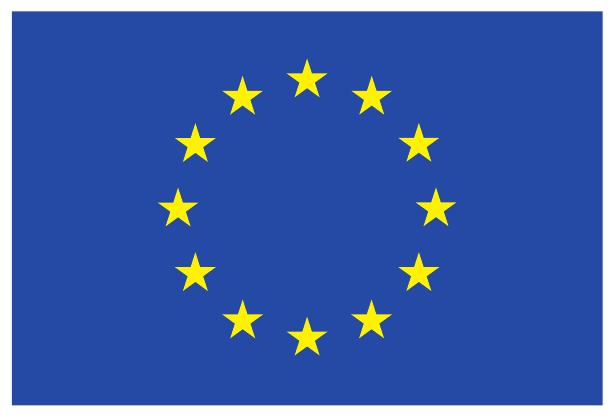} European Union’s Horizon 2020 COFUND programme (LogiCS@TUWien, No. 101034440)}
\and Christopher-Lloyd Simon\affiliationmark{6}%
}

\title{Computing degreewidth of digraphs is hard}
\affiliation{
  DIENS, \'Ecole normale sup\'erieure, CNRS, PSL University, Paris, France\\
  Univ Lyon, CNRS, INSA Lyon, UCBL, Centrale Lyon, Univ Lyon 2, France\\
  Department of Mathematics and Mathematical Statistics, Ume\aa ~University, Sweden\\
  Computer Science Department, Universit\'e libre de Bruxelles, Brussels, Belgium\\
  Algorithm and Complexity Group, TU Wien, Vienna, Austria\\
  The Pennsylvania State University, State College, Pennsylvania, USA%
}

\newcommand\mainthmdetails{}  
\newcommand\bigO[1]{\mathcal{O}(#1)}  

\newcommand\oDelta{\overrightarrow\Delta}

\newcommand\fvn{\mathrm{fvn}}
\newcommand\cost{\mathrm{cost}}
\newcommand\OLA{\mathrm{OLA}}
\newcommand\diOLA{\mathrm{diOLA}}
\newcommand\bw{\mathrm{bw}}
\newcommand\cw{\mathrm{cw}}
\newcommand{\dmax}{d_{\max}}
\newcommand{\dmin}{d_{\min}}
\newcommand{\Dmax}{\Delta_{\max}}
\newcommand{\Dmin}{\Delta_{\min}}
\DeclareMathOperator\dig{dig}
\let\emptyset\varnothing  
\newcommand\card[1]{\left|#1\right|}  

\newcommand\DecisionProblem[1]{\textnormal{\textsc{#1}}\xspace}
\newcommand\DegWidth{\DecisionProblem{Degreewidth}}
\newcommand\dicolourability{\DecisionProblem{Dicolourability}}
\newcommand\FASexp{\DecisionProblem{Feedback Arc Set}}
\newcommand\FAS{\DecisionProblem{FAS}}
\newcommand\cutwidth{\DecisionProblem{Cutwidth}}

\newcommand{\clause}{K}
\newcommand{\cycle}{C}
\newcommand{\lit}{\lambda}

\begin{document}
\publicationdata{vol. 28:2}{2026}{17}{10.46298/dmtcs.14299}{2024-09-18; 2024-09-18; 2026-02-17}{2026-02-25}

\maketitle
\begin{abstract}
Given a digraph, an ordering of its vertices defines a \emph{backedge graph}, namely the undirected graph whose edges correspond to the arcs pointing backwards with respect to the order.
The \emph{degreewidth} of a digraph is the minimum over all ordering of the maximum degree of the backedge graph. We answer an open question by Keeney and Lokshtanov~[WG 2024],
proving that it is {\sf NP}-hard to determine whether an oriented graph has degreewidth at most $1$, which settles the last open case for oriented graphs.
We complement this result with a general discussion on parameters defined using backedge graphs and their relations to classical parameters.
\end{abstract}

\section{Introduction}\label{sec:introduction}
In this work, we consider digraphs $D$ which, unless specified otherwise, will always be \emph{simple}, namely they have no loops and for distinct vertices $u, v\in V(D)$ there is at most one arc from $u$ to $v$ (but $D$ may contain an arc from $u$ to $v$ and an arc from $v$ to $u$).
A \emph{multi-digraph} may have loops and multiple arcs sharing the same source and target.
An \emph{edge} always refers to an undirected graph whereas an \emph{arc} always refers to a digraph.

An \emph{oriented graph} (respectively \emph{tournament}, \emph{semi-complete digraph}) is a digraph in which there is at most (respectively exactly, at least) one arc between each pair of vertices. 
A \emph{symmetric graph} is a directed graph obtained by replacing every edge of an undirected graph by a \emph{digon} (a pair of arcs of opposite orientation). We write $\overleftrightarrow G$ the symmetric graph obtained from $G$.
A \emph{complete symmetric digraph} is a symmetric graph obtained from a complete graph.

Given a digraph $D$ and a total ordering $\prec$ on $V(D)$, we define the associated \emph{backedge graph} as the undirected graph $D^{\prec}$ on the vertex set $V(D)$ whose edges correspond to the arcs $uv\in A(D)$ such that $v\prec u$.

A classical parameter of an undirected graph $G$ is the maximum degree of its vertices, denoted $\Delta(G)$. 
For a digraph $D$, its \emph{degreewidth} $\oDelta(D)$ has recently been introduced in~\cite{davot2023degreewidth} as:
\begin{equation*}
\oDelta(D) = \min \, \left\{ \Delta(D^{\prec}) \mid \mbox{$\prec$ is a total ordering of $V(D)$} \right\}
\end{equation*}

\begin{definition}[$k$-\DegWidth]
For an integer $k\in \mathbb{N}$, we define the problem $k$-\DegWidth as follows:
\begin{itemize}[noitemsep, align=left]
    \item[Input:] a digraph $D = (V, A)$.
    \item[Question:] is $\oDelta(D)$ at most $k$? 
\end{itemize}
\end{definition}

Observe that $0$-\DegWidth is the problem of deciding whether a given digraph is an acyclic digraph, 
and it is clearly in $\P$ since acyclic digraphs are recognisable in linear time using for example a depth first search.

Let us briefly survey the recent work on the complexity of $k$-\DegWidth.

First, for tournaments with $n$ vertices, it is proved in~\cite{davot2023degreewidth} that it is \NP-hard to compute $\DegWidth$, that sorting by in-degree gives a $3$-approximation, and that one can decide if an $n$-vertex tournament has degreewidth $\le 1$ in $\bigO{n^3}$ steps. 

Then, for semi-complete digraphs $D$ with $n$ vertices, it is proved in~\cite{keeney2024degreewidth} that the degreewidth can be computed in $\oDelta(D)^{\bigO{\oDelta(D)}}n + \bigO{n^2}$, and that \FASexp and \cutwidth are $\FPT$ parameterized by the degreewidth on semi-complete digraphs. 

Finally, it is proved in~\cite{keeney2024degreewidth} that deciding if a digraph has degreewidth $\le 2$ is \NP-complete and the complexity of deciding if a digraph has degreewidth $1$ is left as an open question. We solve this question by proving that it is \NP-complete. 

\begin{restatable}[$1$-\DegWidth is \NP-complete]{theorem}{degreewidthNPc}
\label{thm:k-Max-Degree is NP-complete for k >= 1}
For every integer $k \geq 1$, the problem $k$-\DegWidth is \NP-complete\mainthmdetails.
\end{restatable}

The proof will appear in \cref{sec:proof-NPC}. 
In \cref{sec:unifying-undirected-to-directed}, we explain that the definition of $\oDelta$ is not \textit{ad hoc}: it follows a general recipe that we argue should be viewed as a standard way to turn an undirected graph parameter into a directed one (via backedge graphs and minimisation over orderings). We show where this construction already appears and how it provides a unified perspective on several directed parameters. The goal of the section is not to develop new results, but to highlight this unifying viewpoint, which seems to be largely overlooked.


\section{Computing \texorpdfstring{$\oDelta$}{degreewidth} is \NP-complete}
\label{sec:proof-NPC}

\subsection{Preliminary lemmas}

For a digraph $D$, a \emph{feedback arc set} (FAS) is a subset of arcs $F\subseteq A(D)$ such that $D - F$ has no directed cycles. 
Given a digraph $D$ and a subset of arcs $F\subseteq A(D)$, we define the \emph{undirected} graph $D[F]$ with the same vertex set $V(D[F])=V(D)$ and such that for all $u,v\in V(D[F])$, we have $uv \in E(D[F])$ if and only if $(u,v) \in F$ or $(v,u) \in F$. 

\begin{lemma}\label{lemma: degreewidth is equivalent to max deg of a FAS}
    Let $k\in \mathbb{N}$ and $D$ be a digraph.
    We have $\oDelta(D) \leq k$ if and only if $D$ has a FAS $F$ such that $D[F]$ has maximum degree $\Delta(D[F]) \leq k$. 
\end{lemma}

\begin{proof}
    Assume first that $\oDelta(D) \leq k$. There exists an ordering $\prec$ such that $D^{\prec}$ has maximum degree $k$. Hence the set of arcs $F \subseteq A(D)$  corresponding to the edges of $D^{\prec}$ is a feedback arc set such that $\Delta(D[F]) \leq k$. 

    Assume now that $D$ has a feedback arc set $F$ such that $\Delta(D[F])\le k$. 
    The digraph $D-F$ is acyclic and thus admits a topological ordering $\prec$. Hence $E(D^{\prec}) \subseteq F$ and thus $\Delta(D^{\prec})\leq k$. 
\end{proof}

A digraph $D'$ is a \emph{$1$-subdivision} of $D$ if $D'$ can be obtained from $D$ by adding a vertex $v_a$ for every arc $a = (x,y)$ of $D$, and by replacing $a$ by the two arcs $(x,v_a)$ and $(v_a,y)$. An undirected graph is \emph{2-degenerate} if all its subgraphs have a vertex of degree at most $2$. The \emph{underlying graph} of a digraph $D$, is the undirected graph with the same vertices as $D$ and whose edges are obtained from the arcs of $D$ by removing their orientations and ignoring parallel edges.

\begin{remark}
    For a digraph $D$, if $D'$ is a $1$-subdivision of $D$, then the underlying graph $G$ of $D'$ is bipartite and $2$-degenerate.
\end{remark}

\subsection{Proof of \cref{thm:k-Max-Degree is NP-complete for k >= 1}}
\renewcommand\mainthmdetails{ (even restricted to $1$-subdivisions of multi-digraphs)}
Let's prove the main result:
\degreewidthNPc*  

Let $k \geq 1$. The problem $k$-\DegWidth is in \NP, as for a digraph $D$, any given vertex ordering $\prec$ achieving $\Delta(D^{\prec}) \leq k$ can be provided as a certificate.

We prove the \NP-hardness of $k$-\DegWidth by reducing from $3$-\SAT. 

\begin{figure}[ht]
\centering
\begin{tikzpicture}[-,vtx/.style={circle,draw,minimum size=.6cm,text width=.6cm,align=center}]

\node[vtx] (std) at (-8,2){$s$};
\node[vtx] (ttd) at (-4,2){$t$};
\node[td] (td) at (-6,2){$T_p$};
\draw[tdarrow] (std) to (td);
\draw[tdarrow] (td) to (ttd);

\node[] at (-2,2){\large $\Leftrightarrow$};
 
\node[vtx] (1) at (0, 2) {$s$};
\node[vtx] (2) at (2, 4) {$v_1$};
\node[vtx] (3) at (2, 3) {$v_2$};
\node      (dots) at (2, 2) {$\vdots$};
\node[vtx] (t) at (2, 1) {$v_{p-1}$};
\node[vtx] (tp1) at (2, 0) {{$v_p$}};
\node[vtx] (tp2) at (4, 2) {{$t$}};

\begin{scope}[->]
\path (1) edge (2) edge (3) edge (t) edge (tp1);
\path (2) edge (tp2);
\path (3) edge (tp2);
\path (t) edge (tp2);
\path (tp1) edge (tp2);
\end{scope}
\end{tikzpicture}
\caption{Transfer digraph.}
\label{fig:transfer digraph}
\end{figure}
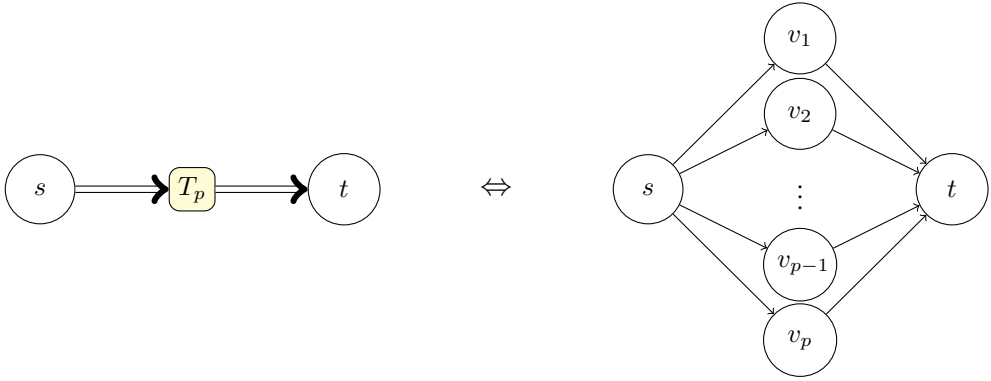

Before starting the reduction, we define for any integer $p$ the {\em transfer digraph of size $p$ from $s$ to $t$} as $T_p = (V^T_p, E^T_p)$  with $ V^T_p=\{s, v_1, \dots, v_{p}, t\}$ and $E^T_p = \{(s, v_i) : 1 \le i \le t\} \cup \{(v_i,t) : 1 \le i \le t\}$. Informally, this digraph consists of $p$ disjoint paths of length $2$ from the source $s$ to the sink $t$, as depicted in \cref{fig:transfer digraph}.

The next claim follows from the  pigeonhole principle. 
\begin{claim}\label{claim delta transfer digraph}
Let $T_{2k} = (V,E)$ be a transfer digraph of size $2k$ from $s$ to $t$ and let $F \subset E$. If $F$ disconnects $s$ from $t$, then $\Delta(D[F]) \ge k$, with equality if and only if $d(s) = d(t) = k$ in $D[F]$. 

Let $T_{2k+1} = (V,E)$ be a transfer digraph of size $2k+1$ from $s$ to $t$ and let $F \subset E$. If $F$ disconnects $s$ from $t$, then $\Delta(D[F]) \ge k+1$.  
\end{claim}

\subsubsection*{Construction of the reduction}

To avoid any confusion, when we write $1 \le i \neq i' \le m$, both the pairs $(i,i')$ and $(i',i)$ are considered. 

On a set of boolean variables $X = \{x_1, \dots x_n\}$, consider a $3$-\SAT\ formula $\varphi = \bigwedge_{1 \leq j \leq m}\clause_j$, where for $1 \le j \le m$ the clause $\clause_j = (\lit_j^1 \vee \lit_j^2 \vee \lit_j^3)$ is a disjunction of literals $\lit_j^i$ for $1 \le i \le 3$.
To the formula $\varphi$, we associate a digraph $D$ as follows.

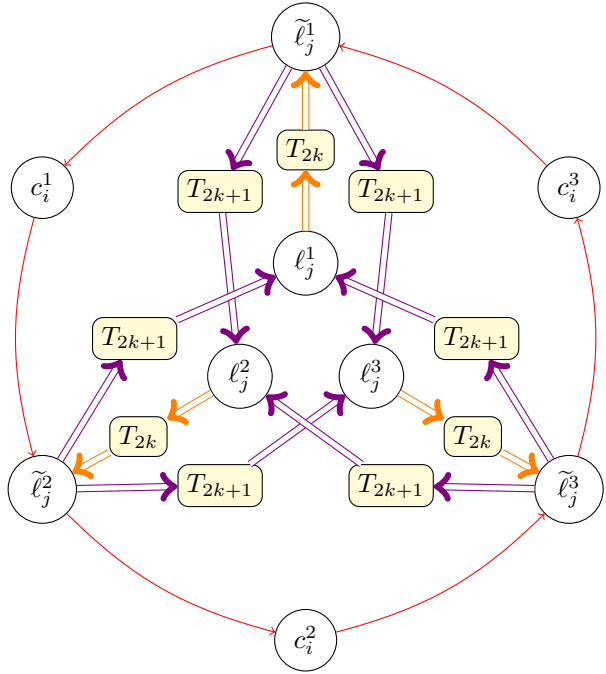
\begin{figure}[ht]
\centering
\begin{tikzpicture}[->, vtx/.style={circle,draw,align=center}]
\node[vtx] (elli1) at (0, 1) {$\ell_j^1$};
\node[vtx] (elli2) at (-0.865, -0.5) {$\ell_j^2$};
\node[vtx] (elli3) at (0.865, -0.5) {$\ell_j^3$};

\node[td] (tdl1) at (0, 2.5) {$T_{2k}$};
\node[td] (tdl2) at (-2.2, -1.3) {$T_{2k}$};
\node[td] (tdl3) at (2.2, -1.3) {$T_{2k}$};

\node[td] (td13) at (1.125, 1.95) {$T_{2k+1}$};
\node[td] (td21) at (-2.25, 0) {$T_{2k+1}$};
\node[td] (td31) at (2.25, 0) {$T_{2k+1}$};
\node[td] (td12) at (-1.125, 1.95) {$T_{2k+1}$};
\node[td] (td23) at (-1.125, -1.95) {$T_{2k+1}$};
\node[td] (td32) at (1.125, -1.95) {$T_{2k+1}$};

\node[vtx] (tildeelli1) at (0, 4) {$\widetilde \ell_j^1$};
\node[vtx] (tildeelli2) at (-3.465, -2) {$\widetilde \ell_j^2$};
\node[vtx] (tildeelli3) at (3.465, -2) {$\widetilde \ell_j^3$};

\node[vtx] (c12) at (-3.465, 2) {{$c_i^{1}$}};
\node[vtx] (c23) at (0, -4) {{$c_i^{2}$}};
\node[vtx] (c31) at (3.465, 2) {{$c_i^{3}$}};

\begin{scope}[every path/.style={draw,color=orange}]

\draw[tdarrow] (elli1) -- (tdl1);
\draw[tdarrow] (tdl1) -- (tildeelli1);

\draw[tdarrow] (elli2) -- (tdl2);
\draw[tdarrow] (tdl2) -- (tildeelli2);

\draw[tdarrow] (elli3) -- (tdl3);
\draw[tdarrow] (tdl3) -- (tildeelli3);

\end{scope}

\begin{scope}[every path/.style={draw,color=violet}]


\draw[tdarrow] (tildeelli1) -- (td13);
\draw[tdarrow] (td13) -- (elli3);


\draw[tdarrow] (tildeelli1) -- (td12);
\draw[tdarrow] (td12) -- (elli2);



\draw[tdarrow] (tildeelli2) -- (td21);
\draw[tdarrow] (td21) -- (elli1);


\draw[tdarrow] (tildeelli2) -- (td23);
\draw[tdarrow] (td23) -- (elli3);



\draw[tdarrow] (tildeelli3) -- (td32);
\draw[tdarrow] (td32) -- (elli2);


\draw[tdarrow] (tildeelli3) -- (td31);
\draw[tdarrow] (td31) -- (elli1);

\end{scope}

\begin{scope}[every path/.style={draw,color=red,bend right=15}]
\path (tildeelli1) edge (c12);
\path (c12)        edge (tildeelli2);
\path (tildeelli2) edge (c23);
\path (c23)        edge (tildeelli3);
\path (tildeelli3) edge (c31);
\path (c31)        edge (tildeelli1);
\end{scope}
\end{tikzpicture}
\caption{Representation of a clause gadget.}
\label{fig:digraph Dj}
\end{figure}

First to every clause $\clause_j = (\lit_j^1 \vee \lit_j^2 \vee \lit_j^3)$ we associate a digraph $D_j$ constructed as follows (see \Cref{fig:digraph Dj}).
\begin{itemize}
    \item For $1 \le i \le 3$, we add the vertices $\ell_j^i$ and $\tilde{\ell}_j^i$ in $D_j$.
    \item For $1 \le i \le 3$, we add a transfer digraph of size $2k$ from $\ell^i_j$ to $\tilde{\ell}^i_j$. 
    \item For $1 \le i \neq i' \le 3$, we add a transfer digraph of size $2k+1$ from $\tilde{\ell}^i_j$ to $\ell_j^{i'}$.
    \item We finally add three vertices $c_j^{1}$, $c_j^{2}$, $c_j^{3}$ and create a directed cycle $\cycle_j = (\tilde{\ell}_j^1, c_j^{1}, \tilde{\ell}_j^2, c_j^{2}, \tilde{\ell}_j^3, c_j^{3}, \tilde{\ell}_j^1)$.
\end{itemize}

We connect the digraphs $D_j$ by adding transfer digraphs as follow:
For each $1 \le j \neq j' \le m$, and each pair $1 \le i,i' \le 3$, such that $\lit_j^i = \neg \lit_{j'}^{i'}$, we add a transfer digraph of size $2k+1$ from $\tilde{\ell}_j^i$ to $\ell_{j'}^{i'}$.  

The resulting digraph $D$ has $\bigO{km^2}$ vertices, so the reduction is polynomial. 
Note that $D$ belongs to the desired class of digraphs: it is a $1$-subdivision of the multi-digraph obtained by replacing every transfer digraph of size $p$ by $p$ parallel arcs and every directed path $(\tilde\ell_j^i, c_j^i, \tilde\ell_j^{i'})$ by an arc $(\tilde\ell_j^i, \tilde\ell_j^{i'})$.

We now prove that $\varphi$ is satisfiable if and only if $D$ has a feedback arc set $F$ with $\Delta(D[F]) \le k$. Dealing with valuation, we use $\top$ ($\bot$ resp.) to specify that a literal is set to 1 (0 resp.). We also identify, the truth value of a literal and its valuation.

\subsubsection*{Intuition behind the reduction}

\begin{enumerate}
    \item During the construction of the digraph $D$, we ensure that every cycle of $D$ either (1) contains all the arcs of $C_j$ (i.e. arcs of the form $(\tilde{\ell},c)$), or (2) goes through a transfer digraph of a literal not chosen to satisfy a clause, which will necessarily be used in the feedback arc set.
    \item The cycle $\cycle_j$ associated to the clause $\clause_j$ guarantees that for every FAS $F$, at least one path joining some $\ell_j^i$ to the associated $\tilde\ell_j^i$ must not intersect $F$, and this index $i$ corresponds to a literal $\lit_j^i$ that must be satisfied (by setting $\nu(x) = \top$ if $\lambda^i_j=x$ and $\nu(x)=\bot$ if $\lambda^i_j=\lnot x$).
    \item The gadget linking literals that are negation of one another ensures that the correspondence in the previous point only provides sound valuations (so that $\nu(\lit) = \top \iff \nu(\lnot\lit) = \bot$ for every literal $\lit$).
\end{enumerate}

\subsubsection*{Proof of the equivalence}

\emph{$(\implies)$}
Suppose first that $\varphi$ is satisfiable. To a valuation  $\nu$ satisfying $\varphi$ we associate a feedback arc sets $F$ such that $\Delta(D[F]) = k$.

We construct the set $F$ as follows:
\begin{itemize}
    \item For each $1 \le j \le m$, we choose one index $i_j$ such that $\nu(\lit_j^{i_j}) = \top$ and we add the arc $(\tilde\ell_j^{i_j}, c_j^{i_j})$ to $F$. Denote by $F_{\top}$ the set of arcs added to $F$ during this step.
    \item For each $1 \le j \le m$, for each $1 \le i \le 3$ with $i \neq i_j$, we add to $F$ $2k$ arcs of the transfer digraph from $\ell_j^{i}$ to $\tilde{\ell}_j^i$ so that these vertices are disconnected.
\end{itemize}

By construction, we have $\Delta(D[F]) = k$. We now prove that $F$ is indeed a feedback arc set. 
Note that for $1 \le j \le m$ and $1 \le i \le 3$, a directed cycle going through a vertex $\ell_j^i$ must also go through $\tilde{\ell}_j^i$. Let $C$ be a directed cycle of $D$.

\begin{itemize}
    \item Suppose first that $C \subset D_j$ for some $1 \le j \le m$. If $C = \cycle_j$ then $A(C) \cap F \neq \emptyset$ since one arc of $\cycle_j$ is in $F_{\top}$. Otherwise, $C$ must contain at least one vertex $\ell_j^i$, and we distinguish two cases.
    \begin{itemize}
        \item If $i \neq i_j$ then $C$ intersects $F$ because all the out-arcs from $\ell_j^i$ are in the transfer digraph from $\ell_j^i$ to $\tilde{\ell}_j^i$, and these have been disconnected  by construction of $F$.
        \item If $i = i_j$, then all out-arcs from $\ell_j^{i_j}$ lead to $\tilde{\ell}_j^{i_j}$. Hence every path from $\tilde\ell_j^{i_j}$ either contains $(\tilde\ell_j^{i_j}, c_j^{i_j})$ which is in $F_\top \subseteq F$, or leads to $\ell_j^i$ for some $i \neq i_j$, and we are back to the previous point. In any case we proved that that $A(C) \cap F \neq \emptyset$.
    \end{itemize}
    
    \item Now suppose that $C$ is not contained in some $D_j$.
    In particular, $C$ must go through a transfer digraph from some $\tilde{\ell}_{j_1}^{i_1}$ to some $\ell_{j_2}^{i_2}$ with $j_1 \neq j_2$.
    
    Recall that by construction of $F$, in $D-F$, for every $j$, there exists a single $i_j$ such that there exists a path between $\ell^{i_j}_j$ and $\tilde\ell_j^{i_j}$, and for that $i_j$, we know that $\nu(\lambda_j^{i_j}) = \top$ and $(\tilde\ell_j^{i_j}, c_j^{i_j}) \in F_\top \subset F$.
    
    By contradiction, suppose that $C \cap F = \emptyset$, i.e. $C$ is a cycle of $D-F$.
    We deduce that $i_2 = i_{j_2}$ and therefore $(\tilde\ell_{j_2}^{i_2}, c_{j_2}^{i_2}) \notin E(C)$ and $\nu(\lambda_{j_2}^{i_2}) = \top$.
    Therefore $C$ must go through the transfer digraph between $\tilde\ell_{j_2}^{i_2}$ and $\ell_{j_3}^{i_3}$.
    Consider the two following cases:
    \begin{itemize}
        \item If $j_2 = j_3$, we know that $i_2 \neq i_3$, and there is no path between $\ell_{j_3}^{i_3}$ and $\tilde\ell_{j_3}^{i_3}$, which is a contradiction.
        \item If $j_2 \neq j_3$, then we know that $\lambda_{j_3}^{i_3} = \lnot \lambda_{j_2}^{i_2}$, hence $\nu(\lambda_{j_3}^{i_3}) = \bot$. In particular $(\tilde\ell_{j_3}^{i_3}, c_{j_3}^{i_3}) \notin F$, therefore all the paths between $\ell_{j_3}^{i_3}$ and $\tilde\ell_{j_3}^{i_3}$ are intersected by $F$, which is a contradiction.
    \end{itemize}
\end{itemize}

\emph{$(\impliedby)$}
Suppose that $D$ has a feedback arc set $F$ with $\Delta(D[F]) \le k$. We construct a valuation $\nu$ that satisfies $\varphi$ as follows. 
For each directed cycle $\cycle_j$, at least one arc needs to be in $F$ and this arc is necessarily incident to exactly one $\tilde{\ell}_j^i$. 
We define $\nu$ such that $\nu(\lit_j^i) = \top$ if an arc from $A(C_j)$ incident to $\tilde{\ell}_j^i$ is in $F$. 
After this operation, if some variables $x_i$ have not been assigned a truth value, then we assign them arbitrarily. 
This construction ensures that for every clause, at least one of its literals is assigned to $\top$. It remains to prove that this valuation is well-defined, namely that there is no literal $\lit$ such that $\nu(\lit) = \top$ and $\nu(\neg \lit) = \top$. 
Suppose by contradiction that there exists some $\lit$ such that $\nu(\lit) = \top$ and $\nu(\neg \lit) = \top$. By definition of $\nu$, there exist two pairs $(j_1, i_1)$ and $(j_2, i_2)$ with $\lit_{j_1}^{i_1} = \neg \lit_{j_2}^{i_2}$ such that an arc containing $\tilde{\ell}_{j_1}^{i_1}$ in $\cycle_{j_1}$ and an arc containing $\tilde{\ell}_{j_2}^{i_2}$ in $\cycle_{j_2}$ are in $F$. 
By construction of the transfer digraph, since one arc of $C_{j_1}$ incident to $\tilde{\ell}_{j_1}^{i_1}$ (respectively an arc of $C_{j_2}$ incident to $\tilde{\ell}_{j_2}^{i_2}$) is in $F$, there exists a path in the transfer digraph from $\ell^{i_1}_{j_1}$ to $\tilde{\ell}^{i_1}_{j_1}$ (respectively from $\ell^{i_2}_{j_2}$ to $\tilde{\ell}^{i_2}_{j_2}$) that does not intersect $F$, otherwise, one of these vertices would have degree at least $k+1$ according to \cref{claim delta transfer digraph}. 
Therefore, since  there is at least one path from $\tilde{\ell}^{i_1}_{j_1}$ to $\ell^{i_2}_{j_2}$ and one path from $\tilde{\ell}^{i_2}_{j_2}$ to $\ell^{i_1}_{j_1}$, there is a directed cycle that does not intersect $F$, contradicting the fact that $F$ is a feedback arc set. Consequently, $\nu$ is well-defined and $\varphi$ is satisfiable.
This concludes the proof of \cref{thm:k-Max-Degree is NP-complete for k >= 1}.\hfill\qed

\section{A Unifying Viewpoint}\label{sec:unifying-undirected-to-directed}

While the definition of $\oDelta$ may seem unusual, it actually follows a general rule which, we believe, should be considered standard for deriving directed graph parameters from undirected graph parameters.

Given an (undirected) graph parameter $\gamma$, we define its directed version as:
\begin{equation}\label{eq:directed_param}
\overrightarrow\gamma(D) = \min \, \left\{ \gamma(D^{\prec}) \mid \mbox{$\prec$ is a total ordering of $V(D)$} \right\}.
\end{equation}
We say that an ordering $\prec$ of a digraph $D$ is a \emph{$\overrightarrow\gamma$-ordering} when $\prec$ witnesses $\overrightarrow\gamma(D)$, that is $\overrightarrow\gamma(D) = \gamma(D^\prec)$.

As for the degreewidth, several digraph parameters are derived by applying \eqref{eq:directed_param} to an undirected graph parameter $\gamma$.
Notably, for certain digraph parameters $\gamma'$, there exists a related undirected parameter $\gamma$ such that $\overrightarrow\gamma$ retrieves $\gamma'$; examples include the dichromatic number and the Feedback Vertex Number, as detailed below. We believe that exploring the perspective given by the backedge graphs can provide new insights into these parameters and potentially uncover novel digraph parameters of interest.

\paragraph{}

In this section, we give an overview of existing digraph parameters defined using~\eqref{eq:directed_param}. 
First, we explain how this construction generalises undirected graph parameters to digraph parameters (via symmetric digraphs) and preserves basic inequalities.
Then, we show that it recovers classical parameters such as the dichromatic number and the feedback vertex number, and we explore  degreewidth a bit further by considering it as a maximum degree concept, rather than merely a parameter for obtaining $\FPT$ algorithms.
Next, we formulate a decision problem that captures the decision problem associated with any digraph parameter defined by~\eqref{eq:directed_param}. Finally, we briefly contrast this framework with layout parameters (optimal linear arrangement, cutwidth and bandwidth), whose standard directed versions are not obtained by applying the same recipe. 


\subsection{Generalisation from graphs to digraphs via symmetric digraphs}


For an undirected graph $G$, the digraph $\overleftrightarrow{G}$ is obtained from $G$ by replacing each edge with a \emph{digon} consisting of two opposite arcs.  
The adjacency relation of $\overleftrightarrow G$ is symmetric, and the map $G\mapsto \overleftrightarrow{G}$ naturally injects graphs into digraphs.
Observe that every ordering $\prec$ of $\overleftrightarrow G$ satisfies ${\overleftrightarrow G}^{\prec} = G$, hence for any graph parameter $\gamma$ we have $\gamma(G) = \overrightarrow\gamma (\overleftrightarrow{G})$. 
Consequently, $\overrightarrow\gamma$ generalises $\gamma$: any property satisfied by $\gamma$ on undirected graphs will also be satisfied by $\overrightarrow\gamma$ in restriction to symmetric digraphs, and one can wonder how this property extends to all other digraphs.  This line of research has been very successful for the dichromatic number (see for example~\cite{M10,HM11,aboulker2023four,Lucas25}), but has received much less attention for other parameters. 
\smallskip 

This generalisation of an undirected graph parameter $\gamma$  to a digraph parameter $\overrightarrow{\gamma}$ also preserves inequalities.
Indeed, if graph parameters $\gamma_1,\gamma_2$ are such that for every graph $G$ we have $\gamma_1(G) \leq \gamma_2(G)$, then for every digraph $D$ and every ordering $\prec$ of $D$ we have $\gamma_1(D^\prec) \leq \gamma_2(D^\prec)$ and thus $\overrightarrow\gamma_1 (D) \leq \overrightarrow\gamma_1(D)$. 
Moreover, since $\gamma\mapsto \overrightarrow{\gamma}$ is defined by a minimum, for a non-decreasing function $f : \mathbb N \to \mathbb R$ we have $\overrightarrow{f\circ \gamma}\le f\circ \overrightarrow{\gamma}$.
Hence monotone inequalities are also preserved: if there exists a non-decreasing function $f : \mathbb R \to \mathbb R$ such that for every graph $G$ we have $\gamma_1(G) \leq f(\gamma_2(G))$, then for every digraph $D$ we have $\overrightarrow\gamma_1(D) \leq \overrightarrow{f\circ \gamma_2}(D) \leq f(\overrightarrow\gamma_2(D))$.
Hence, the hierarchy among undirected parameters is preserved. For example, any tractability ($\FPT$ or $\XP$) result for $\gamma_2$ ensures a similar tractability result for $\gamma_1$, and any hardness ($\W[k]$ or para-$\NP$-hardness) result for $\gamma_1$ implies a similar hardness result for $\gamma_2$.

Finally, note that $\overrightarrow\gamma$ is invariant under reversing the direction of every arc.
For every undirected graph parameter $\gamma$ and digraph $D$ we have $\overrightarrow\gamma(D) = \overrightarrow\gamma(D^{R})$, where $D^{R}$ is obtained from $D$ by reversing the direction of every arc.
Indeed, if $\prec$ is a $\overrightarrow\gamma$-ordering of $D$ then the opposite ordering $\prec'$ defined by $v \prec' w \iff w \prec v$ satisfies $D^\prec=(D^{R})^{\prec'}$ hence $\overrightarrow\gamma(D) \ge \overrightarrow\gamma(D^{R})$ and the other inequality follows by symmetry. 


\subsection{Digraph parameters}



\subsubsection{Dichromatic Number and Directed Clique Number}


We now explain why the directed generalisation $\overrightarrow\chi$ of $\chi$ (obtained by applying \eqref{eq:directed_param} to the chromatic number $\chi$) corresponds to the \emph{dichromatic number} $\chi_a$~\cite{neumann1982dichromatic}, traditionally defined for a digraph $D$ as the minimum $k$ such that $D$ can be partitioned  into $k$ sets, each inducing an acyclic subdigraph.

On the one hand, for every ordering $\prec$, a stable set of $D^{\prec}$ induces an acyclic subdigraph in $D$, so a colouring of $D^{\prec}$ is a dicolouring of $D$ and we have $\chi_a(D) \leq \chi(D^{\prec})$ for every $\prec$. 
On the other hand, it is enough to provide an ordering $\prec$ such that $\chi(D^{\prec}) \leq \chi_a(D)$. 
This can be done by taking an ordering built from a dicolouring,
choosing an ordering on the colour classes, and ordering the vertices of each colour class with a topological order.

This point of view on the dichromatic number has been used several times, for example in~\cite{hero,nguyen2025some,AHP2024}, but we believe it deserves more attention.

A recent development in this direction is the study of the \emph{directed clique number $\overrightarrow\omega$} defined by applying  \eqref{eq:directed_param} to the clique number $\omega$. It was first introduced in~\cite{AACLclique} (see also~\cite{nguyen2025some,A2024}). This provides a notion of clique number for digraphs leading to a notion of $\overrightarrow\chi$-bounded classes of digraphs, which we believe deserve further exploration.

\subsubsection{Directed Vertex Cover is Feedback Vertex Set}

The \emph{vertex cover number} $\tau(G)$ of an undirected graph $G$ is the minimum size of a vertex set $S$ such that  $G - S$ is edgeless.
Its directed version is therefore:
\[
\overrightarrow{\tau}(D) = \min \, \big\{ {\tau}(D^{\prec}) \mid \mbox{$\prec$ is a total ordering of $V(T)$} \big\}
\]
The Feedback Vertex number $\fvn(D)$ of a digraph $D$ is the size of a minimum set of vertices $S$ such that $D - S$ is acyclic. 

Interestingly, $\fvn$ coincides with $\overrightarrow{\tau}$, as we show now.
Let $D$ be a digraph and $S$ be a feedback vertex set of $D$. Then $D-S$ is acyclic, so it has a topological ordering $\prec'$ and thus $(D-S)^{\prec'}$ is edgeless. Extending $\prec'$ to an ordering $\prec$ of $V(D)$ by placing $S$ first yields that $S$ is a vertex cover of $D^{\prec}$, hence $\overrightarrow{\tau}(D)\le |S|=\fvn(D)$. 
Conversely, let $\prec$ be an $\overrightarrow{\tau}(D)$-ordering and let $S$ be a minimum vertex cover of $D^{\prec}$. Since $D^{\prec}-S$ is edgeless, $D-S$ is acyclic, so $\fvn(D)\le |S|=\overrightarrow{\tau}(D)$.

\subsubsection{Back to the Degreewidth}

The maximum degree of a graph has several directed analogues. 
Let $v$ be a vertex of a digraph $D$. 
Recall that for $v \in V(D)$, $d^+(v)$ (resp. $d^-(v)$) denotes the outdegree (resp. indegree) of $v$.
We define the \emph{maxdegree} of $v$ as $\dmax(v)= \max\{d^+(v), d^-(v)\}$ and the \emph{mindegree} of $v$ as $\dmin(v)= \min\{d^+(v), d^-(v)\}$.
We then define the corresponding maximum degrees of $D$: $\Dmax(D)= \max\{\dmax(v)\mid v \in V(D)\}$ and $\Dmin(G)= \max\{\dmin(v)\mid v \in V(D)\}$. The degreewidth $\oDelta$ can be seen as  a third notion of maximum degree. With this point of view, the degreewidth  can be studied by drawing inspiration from results involving the maximum degree in the undirected world.
In this section, we first compare the three notions, and then have a look to a potential generalisation of Brooks' theorem, that relates the maximum degree and the chromatic number of undirected graphs. 

Note that $\oDelta, \Dmax, \Dmin$  coincide on symmetric digraphs. 

\begin{lemma}\label{lem:oriented Delta <= Delta min}
For every digraph $D$, we have $\oDelta(D) \leq \Dmin(D) \leq \Dmax(D)$.
\end{lemma}

\begin{proof}
Let $D$ be a digraph. It is clear that $\Dmin(D) \leq \Dmax(D)$.
Let $\prec$ be an ordering of $V(D)$ that minimises the number of edges of $D^{\prec}$ and assume by contradiction that there exists $u \in V(D)$ such that $d_{D^\prec}(u) > \Dmin(D)$. 
Let $\alpha$ (respectively $\beta$) be the number vertices $x$ such that $x \prec u$ and $(u, x)$ is an arc of $D$ (respectively $u \prec x$ and $(x, u)$ is an arc of $D$).
So $d_{D^\prec}(u) = \alpha + \beta \geq \Dmin(D)+1$. 
Assume first that $d^-(u) = d_{\min}(u)$. Define the ordering $\prec^*$ from $\prec$ by keeping $\prec$ on $V(D) \setminus \{u\}$ untouched and placing $u$ at the beginning. Then
\[
|E(D^{\prec^*})| = |E(D^{\prec})| - \alpha +( d^-(u) - \beta)
\]
Since $d^-(u) \leq \Dmin(D) < \alpha + \beta$, we have $|E(D^{\prec^*})| < |E(D^{\prec})|$, contradicting the choice of $\prec$. 
The case $d^+(u) = d_{\min}(u)$ is treated similarly by placing $u$ at the end of the ordering. 
\end{proof}

\begin{corollary}\label{cor:oriented Delta = Delta min = Delta max for k-regular digraphs}
If $D$ is a $k$-regular digraph, that is $d^-(u) = d^+(u) = k$ for every vertex $u$, then $\oDelta(D) = \Dmin(D) = \Dmax(D) = k$.
\end{corollary}

\begin{proof}
Let $D$ be a $k$-regular digraph. 
It is clear that $\Dmin(D) = \Dmax(D) = k$.
By \cref{lem:oriented Delta <= Delta min}, we know that $\oDelta(D) \leq k$. Moreover, for every ordering $\prec$, the smallest vertex $u$ satisfies $d^-(u) = d_{D^{\prec}}(u) = k$, so $\oDelta(D) \geq k$. 
\end{proof}

Let $D$ be a digraph, and $\prec$ a $\oDelta$-ordering of $D$. By greedily colouring $D^{\prec}$, we get that $\overrightarrow\chi(D) \leq \oDelta(D) + 1$. Hence, for every digraph $D$, we have:

\[\overrightarrow\chi(D) \leq \oDelta(D) +1 \leq \Dmin(D) +1 \leq \Dmax(D) + 1.\]

For $k\in \mathbb{N}$, the $k$-\dicolourability decision problem takes as input a digraph $D$ and decides if it is $k$-dicolourable. 
It is known to be \NP-complete as soon as $k \geq 2$, see~\cite{bokal2004circular}.

Brooks' theorem~\cite{B41} is a fundamental result in graph colouring, stating that every connected graph $G$ satisfies $\chi(G) \leq \Delta(G) +1$ with equality if and only if $G$ is an odd cycle or a complete graph. 
A directed version of Brooks' theorem for $\Dmax$, giving a full characterisation of digraphs $D$ satisfying $\overrightarrow\chi(D) = \Dmax(D) + 1$, has been proved in
~\cite{M10} and a polynomial time algorithm can be deduced to decide if a given digraph $D$ satisfies $\overrightarrow\chi(D) = \Dmax(D) + 1$. In~\cite{aboulker2023four}, it is proved that it is \NP-complete to decide if a given digraph $D$ satisfies $\overrightarrow\chi(D) = \Dmin(D) + 1$.

\begin{theorem}[\cite{aboulker2023four}] \label{thm:Delta_min}
    For $k \geq 2$, $k$-\dicolourability is \NP-complete even when restricted to digraphs $D$ with $\Dmin(D) = k$.
\end{theorem}

Using a very similar proof, we now show that the same holds after replacing $\Dmin$ by $\oDelta$, hence a directed version of Brooks' theorem using $\oDelta$ is unfortunately very unlikely. 

Given a digraph $D$ and a vertex $v$ of $D$, we denote by $\dig_D(v)$ the number of digons incident to $v$ in $D$:
\[\dig_D(v) = \card {\left\{w \in V(D) \mid (v,w), (w,v) \in A(D)\right\}}.\]
If the digraph $D$ is clear, we omit the subscript and write $\dig(v)$.

\begin{lemma}\label{lem:oriented Delta >= max dig(v)}
Every digraph $D$ satisfies $\oDelta(D) \geq \max \; \{\dig(v) \mid v \in V(D)\}$.
\end{lemma}

\begin{proof}
For every $v \in V(D)$, let $w_1, \ldots, w_{\dig(v)}$ be the ``digon-neighbours'' of $v$.
Now fix any ordering $\prec$ of $V(D)$ and observe that for every $1 \leq i \leq \dig(v)$ either $w_i \prec v$, in which case $(v, w_i) \in E(D^{\prec})$, or $v \prec w_i$ in which case  $(w_i, v) \in E(D^{\prec})$. Hence $d_{D^\prec}(v) \geq \dig_D(v)$ which proves the result. 
\end{proof}

\begin{theorem}
\label{thm:k-dicolourabilityy} 
For $k \geq 2$, $k$-\dicolourability is \NP-complete even when restricted to digraphs $D$ with $\oDelta(D) = k$.
\end{theorem}

\begin{proof} 
Let $k \geq 2$. 
In~\cite{aboulker2023four} the proof of \cref{thm:Delta_min} associates to any digraph $D$ a digraph $D'$ satisfying $\Dmin(D') = k$ and $\overrightarrow\chi(D) = k$ if and only $\overrightarrow\chi(D') = k$. We prove that the constructed $D'$ also satisfies $\oDelta(D) = k$, implying the theorem. 
We start, for sake of completeness, by giving the construction of $D'$. 

Start with a given digraph $D$, and construct the digraph $D'$ as follows. 
For each vertex $v$, $D'$ has $k+1$ vertices $v^-, v^+, v_1, \dots v_{k-1}$ (so $D'$ has $|V(D)|(k+1)$ vertices). 
For each vertex $v$ of $D$, create arcs so that $D'[\{v^-, v_1, \dots, v_{k-1}\}]$ and $D'[\{v^+, v_1, \dots, v_{k-1}\}]$ are complete symmetric digraphs, and then add the arc $(v^-, v^+)$.
Finally, for every arc $uv$ of $D$, add the arc $(u^+, v^-)$ to $D'$.

As proved in~\cite{aboulker2023four}, we have $\Dmin(D') = k$. Indeed, for every $v \in V(D)$ and $1 \leq i < k$, we have $\dig_{D'}(v_i) = d^+(v_i) = d^-(v_i) = k$, therefore $\Dmin(D') \geq k$. 
Moreover, for every $v \in V(D)$ we have $d^+(v^-) = d^-(v^+) = k$, so $\dmin(v^+) \leq k$ and $\dmin(v^-) \leq k$. Therefore  $\Dmin(D') = k$.

We now prove that $\oDelta(D) = k$.  Combining \Cref{lem:oriented Delta <= Delta min} and \Cref{lem:oriented Delta >= max dig(v)}, one gets:
\[k \leq \max_{v \in V(D')}\dig(v) \leq \oDelta(D') \leq \Dmin(D') = k.\]

Finally, the proof concludes the same way as for \cref{thm:Delta_min} since $D$ is $k$-dicolourable if and only if $D'$ is.
\end{proof}

Even though we cannot hope for a complete characterisation of digraphs $D$ satisfying  $\overrightarrow\chi(D) = \oDelta(D) +1$, we may still show the following partial result.

\begin{observation}
Let $D$ be a digraph with $\oDelta(D) = k$. If $\overrightarrow\chi(D) = k+1$, then for every $\oDelta$-ordering $\prec$ of $D$, the graph $D^{\prec}$ has a connected component equal to an odd cycle if $k=2$ or to the clique $K_{k+1}$ if $k \neq 2$.  
\end{observation}

\begin{proof}
Let $\prec$ be a $\oDelta$-ordering. If no connected component of $D^{\prec}$  is  one of Brooks' exception, then $\chi(D^{\prec}) \leq \Delta(D^{\prec}) = \oDelta(D)$, and the result follows from $\overrightarrow\chi(D) \leq \chi(D^{\prec})$. 
\end{proof}


\subsection{The decision problems associated with a parameter}


The backedge graph perspective provided by \eqref{eq:directed_param} offers a unified approach to presenting the decision problem related to a digraph parameter $\overrightarrow\gamma$, that is deciding if a given digraph $D$ satisfies $\overrightarrow\gamma(D) \leq k$. 
Given a class of (undirected) graphs $\mathcal C$, we say that a FAS is a $(\mathcal C, \mathcal D)$-FAS if the graph induced by the edges of $F$ (forgetting their orientation) belongs to $\mathcal C$. 
Given a class of
(undirected) graphs $\mathcal C$, we say that a \FAS is a $\mathcal C$-\FAS if the graph induced by the edges of $F$ (forgetting their orientation) belongs to $\mathcal C$. Given a class of undirected graphs $\mathcal C$ and a class of digraphs $\mathcal D$, the \emph{$(\mathcal C, \mathcal D)$-\FAS Problem} is the problem of deciding whether a digraph $D \in \mathcal D$ has a $\mathcal C$-FAS.

In this formalism, the main result of this paper proves that, for every $k \geq 1$, the $(\mathcal C, \mathcal D)$-FAS problem in oriented graphs (more precisely in $1$-subdivision of oriented graphs) is \NP-hard when $\mathcal C$ is the class of graphs with maximum degree $k$. 

More generally, for a graph parameter $\gamma$ which is monotone (in the sense $\gamma(H) \leq \gamma(G)$ for any subgraph $H$ of $G$), deciding whether a digraph $D$ satisfies $\overrightarrow\gamma(D) \leq k$ is the same as the $(\mathcal C, \mathcal D)$-FAS problem in $\mathcal D$ where $\mathcal C$ is the class of graphs $G$ satisfying $\gamma(G) \leq k$ and $\mathcal D$ is the class of all digraphs. This problem is often of particular interest when $\mathcal D$  is taken to be the class of all tournaments. For more on this subject we refer to~\cite{AAL2024} (in particular the last section).

\subsection{Optimal Linear Arrangement, Cutwidth and Bandwidth}
\label{subsec:OLA-cutwidth-bandwidth}



Given an ordered undirected graph $(G, \prec)$, the \emph{length} of an edge $uv$ such that $u \prec v$ is defined as $1$ plus the number of vertices $x$ such that $u \prec x \prec v$.
One can define three undirected parameters by minimising certain functions of the lengths over all ordering of the graph: Optimal Linear Arrangement, cutwidth, and bandwidth.
Each of these parameters has a directed version defined similarly, but for a given digraph $D$ and an ordering $\prec$, only the edges of $D^{\prec}$ are considered.
These directed versions differ from those obtained by applying \eqref{eq:directed_param} to the undirected parameters.
The sole purpose of this section is to highlight this difference.
The authors did not explore the relationship between these parameters.

Let us provide a detailed explanation for the Optimal Linear Arrangement. 
For an undirected graph $G$ with an ordering $\prec$, we denote by $\cost(G, \prec)$ the sum of the lengths of its edges.
An \emph{Optimal Linear Arrangement} of $G$ is an ordering $\prec$ that minimises $\cost(G, \prec)$, and its value is denoted $\OLA(G)$, thus:
\begin{equation*}
\OLA(G) = \min\;\{\cost(G, \prec) \mid \mbox{$\prec$ is a total ordering of $V(G)$}\}.
\end{equation*}

The directed version $\diOLA$ of $\OLA$ is defined similarly (see~\cite{CCDG82,FP19}), but only the length of backward arcs are summed. More precisely, 
\begin{equation*}
\diOLA(D) = \min \; \{\cost(D^{\prec}, \prec) \mid \mbox{$\prec$ is a total ordering of $V(D)$} \}.
\end{equation*}
This is not the same as $\overrightarrow{\OLA}$ defined according to the formula \eqref{eq:directed_param} by:
\begin{align*}
    \overrightarrow{\OLA}(D) 
    & = \min \; \{\OLA(D^{\prec}) \mid \mbox{$\prec$ is a total ordering of $V(D)$}\} \\
    & = \min \; \{\cost(D^{\prec}, \prec') \mid \mbox{$\prec, \prec'$ are total orderings of $V(D)$}\}.
\end{align*}
Hence, in $\overrightarrow{\OLA}$, we can choose the Optimal Linear Arrangement over each backedge graph $D^{\prec}$, while in $\diOLA$ we only consider the cost of $(D^{\prec}, \prec)$.
Thus, for every digraph $D$, we have $\overrightarrow{\OLA}(D) \leq \diOLA(D)$.

The \emph{bandwidth} of an ordered graph $(G, \prec)$ is the maximum length of an edge, and the bandwidth of a graph $G$ is  $\bw(G) = \min\{\bw(G, \prec) \mid \mbox{$\prec$ is a total ordering of $V(G)$}\}$. 
The \emph{cutwidth} $\cw(G)$ of an undirected graph $G$ is the smallest integer $k$ with the following property: there is an ordering of the vertices of the graph, such that every cut obtained by partitioning the vertices into earlier and later subsets of the ordering is crossed by at most $k$ edges.
As with $\OLA$, the directed bandwidth~\cite{JKLSS19} and directed cutwidth~\cite{CFS12} of a digraph $D$ are defined similarly to their undirected counterparts, but for each ordering $\prec$, only the edges of $D^{\prec}$ are considered. Therefore, as with $\OLA$, applying \eqref{eq:directed_param} to these parameters does not yield the same results.

\acknowledgements\label{sec:ack}
We are grateful to the two anonymous referees for their helpful comments on the first version of the manuscript.

\nocite{*}
\bibliographystyle{abbrvnat}
\bibliography{references}

@article{neumann1982dichromatic,
  title={The dichromatic number of a digraph},
  author={Neumann-Lara, Victor},
  journal={Journal of Combinatorial Theory, Series B},
  volume={33},
  number={3},
  pages={265--270},
  year={1982},
  publisher={Elsevier}
}

@article{bokal2004circular,
  title={The circular chromatic number of a digraph},
  author={Bokal, Drago and Fijavz, Gasper and Juvan, Martin and Kayll, P Mark and Mohar, Bojan},
  journal={Journal of Graph Theory},
  volume={46},
  number={3},
  pages={227--240},
  year={2004},
  publisher={Wiley Online Library}
}

@article{aboulker2023four,
  title={Four proofs of the directed {B}rooks' {T}heorem},
  author={Aboulker, Pierre and Aubian, Guillaume},
  journal={Discrete Mathematics},
  volume={346},
  number={11},
  pages={113193},
  year={2023},
  publisher={Elsevier}
}

@inproceedings{davot2023degreewidth,
  title={{D}egreewidth: {A} {N}ew Parameter for {S}olving {P}roblems on {T}ournaments},
  author={Davot, Tom and Isenmann, Lucas and Roy, Sanjukta and Thiebaut, Jocelyn},
  booktitle={International Workshop on Graph-Theoretic Concepts in Computer Science},
  pages={246--260},
  year={2023},
  organization={Springer}
}

@inproceedings{keeney2024degreewidth,
  title={{D}egreewidth on semi-complete digraphs},
  author={Keeney, Ryan and Lokshtanov, Daniel},
  booktitle={International Workshop on Graph-Theoretic Concepts in Computer Science},
  pages={312--326},
  year={2024},
  organization={Springer}
}

@article{nguyen2025some,
  title={Some results and problems on tournament structure},
  author={Nguyen, Tung and Scott, Alex and Seymour, Paul},
  journal={Journal of Combinatorial Theory, Series B},
  volume={173},
  pages={146--183},
  year={2025},
  publisher={Elsevier}
}

@article{B41,
	title        = {On colouring the nodes of a network},
	author       = {R. L. Brooks},
	year         = 1941,
	journal      = {Math. Proc. Cambridge Philos. Soc.},
	volume       = 37,
	pages        = {194--197}
}

@article{CCDG82,
    author = {Chinn, P. Z. and Chvátalová, J. and Dewdney, A. K. and Gibbs, N. E.},
    title = {The bandwidth problem for graphs and matrices—a survey},
    journal = {Journal of Graph Theory},
    volume = {6},
    number = {3},
    pages = {223-254},
    doi = {10.1002/jgt.3190060302},
    url = {https://onlinelibrary.wiley.com/doi/abs/10.1002/jgt.3190060302},
    eprint = {https://onlinelibrary.wiley.com/doi/pdf/10.1002/jgt.3190060302},
    year = {1982}
}

@article{M10,
    title        = {Eigenvalues and colourings of digraphs},
    author       = {B. Mohar},
    year         = 2010,
    journal      = {Linear Algebra and its Applications},
    volume       = 432,
    number       = 9,
    pages        = {2273--2277}
}

@article{HM11,
  title={Gallai's theorem for list coloring of digraphs},
  author={Harutyunyan, Ararat and Mohar, Bojan},
  journal={SIAM Journal on Discrete Mathematics},
  volume={25},
  number={1},
  pages={170--180},
  year={2011},
  publisher={SIAM}
}

@article{CFS12,
  title={Tournament immersion and cutwidth},
  author={Chudnovsky, Maria and Fradkin, Alexandra and Seymour, Paul},
  journal={Journal of Combinatorial Theory, Series B},
  volume={102},
  number={1},
  pages={93--101},
  year={2012},
  publisher={Elsevier}
}

@article{hero,
author = {Berger, Eli and Choromanski, Krzysztof and Chudnovsky, Maria and Fox, Jacob and Loebl, Martin and Scott, Alex and Seymour, Paul and Thomassé, Stéphan},
year = {2013},
month = {01},
pages = {1–20},
title = {Tournaments and colouring},
volume = {103},
journal = {Journal of Combinatorial Theory, Series B},
doi = {10.1016/j.jctb.2012.08.003}
}

@article{FP19,
title = {On width measures and topological problems on semi-complete digraphs},
journal = {Journal of Combinatorial Theory, Series B},
volume = {138},
pages = {78-165},
year = {2019},
issn = {0095-8956},
doi = {10.1016/j.jctb.2019.01.006},
url = {https://www.sciencedirect.com/science/article/pii/S0095895619300073},
author = {Fedor V. Fomin and Michał Pilipczuk}
}

@inproceedings{JKLSS19,
  author       = {Pallavi Jain and
                  Lawqueen Kanesh and
                  William Lochet and
                  Saket Saurabh and
                  Roohani Sharma},
  editor       = {Arkadev Chattopadhyay and
                  Paul Gastin},
  title        = {Exact and Approximate Digraph Bandwidth},
  booktitle    = {39th {IARCS} Annual Conference on Foundations of Software Technology  and Theoretical Computer Science, {FSTTCS} 2019, December 11-13, 2019,   Bombay, India},
  series       = {LIPIcs},
  volume       = {150},
  pages        = {18:1--18:15},
  publisher    = {Schloss Dagstuhl - Leibniz-Zentrum f{\"{u}}r Informatik},
  year         = {2019},
  url          = {10.4230/LIPIcs.FSTTCS.2019.18},
  doi          = {10.4230/LIPICS.FSTTCS.2019.18},
  bibsource    = {dblp computer science bibliography, https://dblp.org}
}

@article{AACLclique,
  title={Clique number of tournaments},
  author={Aboulker, Pierre and Aubian, Guillaume and Charbit, Pierre and Lopes, Raul},
  journal={arXiv preprint arXiv:2310.04265},
  year={2023}
}

@article{A2024,
  title={Computing the clique number of tournaments},
  author={Aubian, Guillaume},
  journal={arXiv preprint arXiv:2401.07776},
  year={2024}
}

@article{AAL2024,
  title={Finding forest-orderings of tournaments is {NP}-complete},
  author={Aboulker, Pierre and Aubian, Guillaume and Lopes, Raul},
  journal={arXiv preprint arXiv:2402.10782},
  year={2024}
}

@article{AHP2024, 
title={Minimum Acyclic Number and Maximum Dichromatic Number of Oriented Triangle-Free Graphs of a Given Order}, volume={32}, 
url={https://www.combinatorics.org/ojs/index.php/eljc/article/view/v32i4p27}, 
DOI={10.37236/12862}, 
number={4}, 
journal={The Electronic Journal of Combinatorics}, 
author={Aboulker, Pierre and Havet, Frédéric and Pirot, François and Schabanel, Juliette}, year={2025}, month={Nov.}, pages={P4.27} }

@inproceedings{Lucas25,
  title={An analogue of Reed’s conjecture for digraphs},
  author={Kawarabayashi, Ken-Ichi and Picasarri-Arrieta, Lucas},
  booktitle={Proceedings of the 2025 Annual ACM-SIAM Symposium on Discrete Algorithms (SODA)},
  pages={3310--3324},
  year={2025},
  organization={SIAM}
}
\label{sec:biblio}

\end{document}